\documentstyle[12pt]{article}
\def\tilde{\widetilde}

\def\r{r_{+}} 
\def\Rp{R_{+}} 
\def\Rpi{R_{+}^{-1}} 
\def\Ra{R_{-}} 
\def\Rai{R_{-}^{-1}} 
\def\Rb{R_{\pm}} 
\def\ra{r_{-}} 
\def\rb{r_{\pm}}

\def\trq{\,{\rm tr_q}\,} 
\def\detq{\,{\rm det_q}\,} 
\def\a{\alpha} 
\def\b{\beta}

\def\l{\lambda} 
\def\G.{Gauss-law constraints. } 
\def\G{Gauss-law constraints }

\def\la{\lambda^{a}} 
\def\lb{\lambda^{b}}

\def\be{\begin{equation}} 
\def\ee{\end{equation}} 
\def\bea{\begin{eqnarray}} 
\def\eea{\end{eqnarray}}

\def\a{\alpha}

\def\Hs{Heisenberg doubles } 
\def\H.{the Heisenberg double. } 
\def\ah.{a Heisenberg double. } 
\def\Hs.{Heisenberg doubles. }

\def\Lg{{\cal L}_g}
\def\lo{{\cal L}_1}
\def\Lo{{\cal L}_1}

\def\Lot{{\cal L}_1 (sl_2)}
\def\ELg{{\cal L}_{g}^{*}}
\def\ELo{{\cal L}_{1}^{*}}
\def\PLg{{\cal PL}_g}
\def\PLo{{\cal PL}_1}

\def\lt{{\cal L}_g(sl_2)}
\def\Lt{{\cal L}_g(sl_2)}
\def\Mg{{\cal M}_g}
\def\PM{{\cal PM}_g}
\def\Mo{{\cal M}_g(sl_2)}
\def\cMo{{\cal Z}(\Mo)}
\def\F{{\cal F}_\I}
\def\I{{\cal I}}
\vspace{5cm}
\title{  \hfill{LMU-TPW 96-11} \\ 
\hfill{q-alg/9603...}\\ 
\vspace{1cm}
The centre of the graph and moduli algebras at roots of 1.}
\vspace{.7cm}
\author{ \mbox{}\\ 
S.A.Frolov\thanks{Alexander von Humboldt fellow} \mbox{} \\ 
\vspace{0.4cm} 
Section Physik, Munich University 
\vspace{-0.5cm} \mbox{} \\ 
Theresienstr.37, 80333 Munich, Germany
\thanks{Permanent address:\ Steklov Mathematical Institute, Vavilov st.42,
 GSP-1, 117966 Moscow, RUSSIA}
\mbox{}}  
\date{}
\begin{document}
\maketitle
\vspace{3.5cm}
\begin{abstract}
The structure of the centres ${\cal Z}(\Lg)$ and ${\cal Z}(\Mg)$ of the graph 
algebra ${\cal L}_g(sl_2)$ and the moduli algebra ${\cal M}_g(sl_2)$ is studied 
at roots of 1. It it shown that ${\cal Z}(\Lg)$ can be endowed with the
structure of the Poisson graph algebra. The elements of $Spec({\cal Z}(\Mg))$ 
are shown to satisfy the defining relation for the holonomies of a flat 
connection along the cycles of a Riemann surface. The irreducible 
representations of the graph algebra are constructed.
\end{abstract} 

\section{Introduction} 

The Poisson structure of the moduli space of flat connections on a Riemann
surface with $g$ handles can be described by means of a quadratic Poisson
algebra, which was introduced by Fock and Rosly \cite{FR} and here will be 
called
the Poisson graph algebra. Let us remind the definition of the algebra
\cite{FR}.

Let $G$ be a matrix algebraic group and $D_g =G^{\times^{2g}}$. An arbitrary
element $d$ of $D_g$ is parametrized by matrices $A_i$ and $B_i$ as 
$d=(A_1 ,B_1 ,...,A_g ,B_g )\in D_g$. Let all of the matrices be in the
fundamental representation of the group $G$. Then the algebra of functions on
$D_g$ is generated by the matrix elements $(A_i )_{mn}$ and  $(B_i )_{mn}$.

{\bf Definition 1.} The Poisson graph algebra $\PLg$ is an algebra of regular
functions on $D_g$ with the following Poisson structure
\bea &&i=1,\cdots ,g \nonumber\\ 
\frac {k}{2\pi}\{ A_{i}^{1},A_{i}^{2}\}&=&
A_{i}^{1}\r A_{i}^{2}- A_{i}^{2}A_{i}^{1}\r -\ra 
A_{i}^{2}A_{i}^{1}+A_{i}^{2}\ra A_{i}^{1} \nonumber\\ 
\frac {k}{2\pi}\{ B_{i}^{1},B_{i}^{2}\}&=&
B_{i}^{1}\r B_{i}^{2}- B_{i}^{2}B_{i}^{1}\r -\ra 
B_{i}^{2}B_{i}^{1}+B_{i}^{2}\ra B_{i}^{1} \nonumber\\ 
\frac {k}{2\pi}\{ A_{i}^{1},B_{i}^{2}\}&=&
A_{i}^{1}\r B_{i}^{2}- B_{i}^{2}A_{i}^{1}\r -\r 
B_{i}^{2}A_{i}^{1}+B_{i}^{2}\ra A_{i}^{1} \nonumber\\ &&i<j \nonumber\\ 
\frac {k}{2\pi}\{ A_{i}^{1},A_{j}^{2}\}&=&A_{i}^{1}\r A_{j}^{2}-
A_{j}^{2}A_{i}^{1}\r -\r A_{j}^{2}A_{i}^{1}+A_{j}^{2}\r A_{i}^{1}
\nonumber\\
\frac {k}{2\pi}\{ A_{i}^{1},B_{j}^{2}\}&=&A_{i}^{1}\r B_{j}^{2}-
B_{j}^{2}A_{i}^{1}\r -\r B_{j}^{2}A_{i}^{1}+B_{j}^{2}\r A_{i}^{1}
\nonumber\\
\frac {k}{2\pi}\{ B_{i}^{1},B_{j}^{2}\}&=&B_{i}^{1}\r B_{j}^{2}-
B_{j}^{2}B_{i}^{1}\r -\r B_{j}^{2}B_{i}^{1}+B_{j}^{2}\r B_{i}^{1}
\nonumber\\
\frac {k}{2\pi}\{ B_{i}^{1},A_{j}^{2}\}&=&B_{i}^{1}\r A_{j}^{2}-
A_{j}^{2}B_{i}^{1}\r -\r A_{j}^{2}B_{i}^{1}+A_{j}^{2}\r B_{i}^{1}
\label{CPL}
\eea
Here $k$ is an arbitrary complex parameter, $\rb$ are classical $r$-matrices 
which satisfy the classical Yang-Baxter equation and the following relations
\be
[r^{12} ,r^{13} ]+[r^{12} ,r^{23} ]+[r^{13} ,r^{23} ]=0
\label{cyb1}
\ee
\be
\ra =-P\r P, \qquad \r -\ra =C
\label{cyb}
\ee
where $P$ is a permutation in the tensor product $V\otimes V$ ($Pa\otimes 
b=b\otimes a$).

\noindent In eq.(\ref{CPL}-\ref{cyb}) we use the standard notations from the 
theory of 
quantum groups \cite{D,FRT}:  for any matrix $A$ acting in some space $V$ 
one can construct two matrices $A^{1}=A\otimes id $ and $A^{2}=id\otimes A $ 
acting in the space $V\otimes V$ and  for any matrix $r=\sum_{a} r_{1}(a)
\otimes r_{2}(a)$ acting in the space $V\otimes V$ one constructs matrices 
$r^{12}=\sum_{a} r_{1}(a)\otimes r_{2}(a)\otimes id$, 
$r^{13}=\sum_{a} r_{1}(a)\otimes id\otimes r_{2}(a)$ and 
$r^{23}=\sum_{a} id\otimes r_{1}(a)\otimes r_{2}(a)$ acting in the space 
$V\otimes V\otimes V$. The matrix $C$ is the tensor Casimir 
operator of the Lie algebra ${\cal G}$ of the group $G$: 
$C=-\eta_{ab}\la\otimes\lb$, $\eta_{ab}$ is the Killing tensor and 
$\la$ form a basis of ${\cal G}$.

Let us now identify the matrices $A_i$ and $B_i$ with holonomies of a flat
connection along the cycles $a_i$ and $b_i$ of a Riemann surface with $g$
handles. Then $A_i$ and $B_i$ should satisfy the following defining relations
\be
M=B_{g}A_{g}^{-1}B_{g}^{-1}A_{g}\cdots B_{1}A_{1}^{-1}B_{1}^{-1}A_1 = 1
\label{R}
\ee
These relations can be regarded as first-class constraints imposed on the
variables of $D_g$. The gauge transformations generated by these constraints
are just the simultaneous conjugations of $A_i$ and $B_i$
\[ A_i\to hA_i h^{-1},\qquad A_i\to hA_i h^{-1}\]

Let us now consider two Poisson subalgebras of $\PLg$. The first subalgebra $\I$
consists of all of the functions vanishing on the constraints surface: 
\[
\I =\{ f\in \PLg : f| _{M=1} =0\}
\]
and the second subalgebra $\F$ is the maximal subalgebra of $\PLg$ such that the
subalgebra $\I$ is a Poisson ideal of $\F$:
\[
\F =\{ f\in \PLg : \{ f,h\}\in \I \quad \forall h\in \I\} .
\]
In particular it is not difficult to check that any function $f$ which is
invariant with respect to simultaneous conjugations of $A_i$ and $B_i$
\[
f(hA_1 h^{-1},...,hB_g h^{-1})=f(A_1 ,...,B_g )
\]
belongs to $\F$.

{\bf Definition 2.} The Poisson algebra of functions on the moduli space of flat
connections on a Riemann surface with $g$ handles or the Poisson moduli algebra
$\PM$ is defined as a quotient of the algebra $\F$ over the ideal $\I$
\[
\PM =\F \big/ \I.
\]
It was shown by Fock and Rosly that the algebra $\PM$ coincides with the
canonical Poisson algebra of functions on the moduli space defined by the
Atiyah-Bott symplectic structure.

Quantization of the Poisson graph algebra leads to an associative algebra which
was introduced in \cite{AGS} (see also \cite{BR}). In the present paper we study
 the structure of the
centre of the quantized graph and moduli algebras for the simplest case of
$SL(2)$ group and, in what follows, present definitions and results only for
this case. Our definition of the moduli algebra differs from the definition
given in \cite{AGS}, where the truncated case was considered, and can be
regarded as the standard one from quantum theory of constraints systems.

The plan of the paper is as follows. In the second section we introduce the 
graph
and moduli algebras. Then we describe an extension of the graph algebra and the 
isomorphism between the extended graph algebra $\ELg$ and the tensor product of
$g$ copies of $\ELo$ \cite{AM,A}. In the third section we study the centre
${\cal Z}(\lo)$ of $\lo$ and prove that ${\cal Z}(\lo)$ is isomorphic to
$\PLo$. In the forth section,
using the isomorphism mentioned in the second section, we generalize the results
obtained in the third section to $\Lg$ and show that  the elements
of $Spec(\cMo)$ satisfy the defining relation (\ref{R}). In the fifth section 
the irreducible
representations of $\Lg$ are constructed. In Conclusion we discuss unsolved
problems.

\section{Graph and moduli algebras}  

{\bf Definition 3.} The graph algebra $\lt$ is an associative algebra with unit
element, generated by matrix elements of $A_i ,B_i \in End \, C^2 \otimes 
\Lt , \quad
i=1,...,g$ and $(A_i )_{11}^{-1}$, $(B_i )_{11}^{-1}$, $M_{11}^{-1}$ which
 are subject to the     
following relations 
\[ i=1,\cdots ,g \]
\[ A_{i}^{1}\Rp A_{i}^{2}\Rpi=\Ra A_{i}^{2}\Rai A_{i}^{1},\quad 
B_{i}^{1}\Rp B_{i}^{2}\Rpi=\Ra B_{i}^{2}\Rai B_{i}^{1}\]
\be
A_{i}^{1}\Rp B_{i}^{2}\Rpi=\Rp B_{i}^{2}\Rai A_{i}^{1} 
\label{ab}
\ee
\bea
&&\qquad\qquad\qquad\qquad i<j \nonumber\\ 
&&A_{i}^{1}\Rp A_{j}^{2}\Rpi=\Rp A_{j}^{2}\Rpi A_{i}^{1},\quad
A_{i}^{1}\Rp B_{j}^{2}\Rpi=\Rp B_{j}^{2}\Rpi A_{i}^{1} \nonumber\\ 
&&B_{i}^{1}\Rp B_{j}^{2}\Rpi=\Rp B_{j}^{2}\Rpi B_{i}^{1},\quad 
B_{i}^{1}\Rp A_{j}^{2}\Rpi=\Rp A_{j}^{2}\Rpi B_{i}^{1} \nonumber\\ 
&&(A_i )_{11}(A_i )_{11}^{-1}=(A_i )_{11}^{-1}(A_i )_{11}=
(B_i )_{11}(B_i )_{11}^{-1}=(B_i )_{11}^{-1}(B_i )_{11}=1\nonumber\\
&&M_{11}^{-1}
(q^{-3g}B_{g}A_{g}^{-1}B_{g}^{-1}A_{g}\cdots B_{1}A_{1}^{-1}B_{1}^{-1}A_1)_{11}
=\nonumber\\
&&(q^{-3g}B_{g}A_{g}^{-1}B_{g}^{-1}A_{g}\cdots B_{1}A_{1}^{-1}B_{1}^{-1}
A_1)_{11}M_{11}^{-1}=1\nonumber\\
&&\detq A_i =(A_i )_{11}(A_i )_{22}-q^2 (A_i )_{21}(A_i )_{12} =1\nonumber\\
&&\detq B_i =(B_i )_{11}(B_i )_{22}-q^2 (B_i )_{21}(B_i )_{12} =1
\label{QR}
\eea
and we denote by 1 the unit element of any algebra throughout the paper.

Here $\Rb$-matrices 
\[
\Rp =q^{-\frac {1}{2}}\left( \begin{array}{cccc} q & 0 &  0 & 0 \\
0 & 1 & q-q^{-1} & 0 \\ 
0 & 0 &  1 & 0 \\ 
0 & 0 & 0 & q \end{array} \right),\qquad 
\Ra =q^{\frac {1}{2}}\left( \begin{array}{cccc} q^{-1} & 0 &  0 & 0 \\
0 & 1 & 0 & 0 \\ 
0 & q^{-1}-q &  1 & 0 \\ 
0 & 0 & 0 & q^{-1} \end{array} \right)
\]
satisfy the quantum Yang-Baxter equation and the following
relations
\[
R^{12}R^{13}R^{23}=R^{23}R^{13}R^{12} 
\]
\[
\Rp =P\Rai P ,\qquad 
\Rb (q) = 1 + \frac {2\pi i}{k} \hbar \rb +O(\hbar^2) \]
\[ q=\exp ( \frac {2\pi i}{p} \hbar ), \qquad p=k +2\hbar
\]
We have introduced the element $M_{11}^{-1}$ because the Poisson structure
(1.1) is degenerate on the surface $M_{11}=0$ and because, as we will see in
the fifth section, if the element $M_{11}$ is invertible then all irreducible
representations of the graph algebra are $p^{3g}$-dimensional.

Using the explicit expressions for $\Rb$ one can write down the commutation
relations (\ref{QR}) in terms of elements $(A_i)_{mn}$ and  $(B_i)_{mn}$. The
corresponding formulas are given in Appendix.

It is well known that for fixed index $i$ the algebra generated by matrix
elements of $A_i$ (or $B_i$) is isomorphic to $U_q (sl_2)$ \cite{RS} and the
standard generators of $U_q (sl_2)$ can be expressed through $(A_i)_{mn}$ as
follows
\[
K=(A_i)_{11},\quad X_+ =-\frac {q^{\frac {1}{2}}} {q-q^{-1}}(A_i)_{21},\quad 
X_- =-\frac {q^{\frac {1}{2}}} {q-q^{-1}}(A_i)_{11}^{-1}(A_i)_{12}
\]
The Casimir element of $U_q (sl_2)$ is equal to
\[
c_i =\trq A_i =q^{-1}(A_i)_{11}+q(A_i)_{22}
\]
However the Casimir elements $c_i$ are not central elements of the graph algebra
and, moreover, for generic values of $q$ the graph algebra has a trivial centre.

It is not difficult to show that for fixed index $i$ the algebra generated by 
$A_i$ and $B_i$ is isomorphic to the quantized algebra of functions on the
Heisenberg double of a Lie group \cite{S1,AF1,S2,AF2}.

Let us introduce monodromies $M_i$
\[
M_i =q^{-3(g-i+1)}B_{g}A_{g}^{-1}B_{g}^{-1}A_{g}\cdots 
B_{i}A_{i}^{-1}B_{i}^{-1}A_i
\]
The monodromies satisfy the following commutation relations
\[ i\leq j \] 
\[ M_{i}^{1}\Rp M_{j}^{2}\Rpi=\Ra M_{j}^{2}\Rai M_{i}^{1},\quad
M_{i}^{1}\Rp A_{j}^{2}\Rpi=\Ra A_{j}^{2}\Rai M_{i}^{1} \] 
\[ M_{i}^{1}\Rp B_{j}^{2}\Rpi=\Ra B_{j}^{2}\Rai M_{i}^{1} \] 
\[ i<j \] 
\[
A_{i}^{1}\Rp M_{j}^{2}\Rpi=\Rp M_{j}^{2}\Rpi A_{i}^{1},\quad 
B_{i}^{1}\Rp M_{j}^{2}\Rp=\Rp M_{j}^{2}\Rpi B_{i}^{1}\]
\be
\detq M_i =1
\label{CRM}
\ee
Using the Gauss decomposition one can represent the monodromies as follows
\[ M_i =m_{-}^{-1}(i)K_i m_+ (i)\]
 where $m_{\pm}(i)$ are upper- and lower-triangular matrices with the
unity on the diagonals and $K_i$ are diagonal matrices. It follows from
relations (\ref{CRM}) that the matrix elements of $K_i$ form a commutative
subalgebra of $\Lg$. 

Let $\ELg$ be an extension of the graph algebra $\Lg$ by means of the
elements $Q_{i}^{\pm 1}=K_{i}^{\pm \frac {1}{2}}$. The relations (\ref{CRM}) and
the commutation relations of $Q_i$ are presented in components in Appendix.
Then the following proposition is a refinement of some results from \cite{AM,A}:

{\bf Proposition 1.} The extended graph algebra $\ELg$ is isomorphic to the
tensor product of $g$ copies of the extended graph algebra $\ELo$. The
isomorphism is given by means of the following formulas
\be
A_i =M_{+}^{-1}(i+1)\bar A_i M_{+}(i+1),\qquad 
B_i =M_{+}^{-1}(i+1)\bar B_i M_{+}(i+1)
\label{IS}
\ee
where
\[ M_+ (i)=Q_i m_+ (i), \qquad M_- (i)=Q_i^{-1}m_{-}(i)\]
\[ M_{\pm}(i)=\bar G_{\pm}(i)\bar G_{\pm}(i+1)\cdots\bar G_{\pm}(g)\] 
\[ \bar G_{i}=q^{-3}\bar B_i\bar A_i^{-1}\bar B_i^{-1}\bar A_i \]

{\bf Remark 1.} The Proposition 1 is valid for arbitrary graph algebra. In the
case of $SL(2)$ group the element $(K_i )_{11}=(M_i )_{11}=(\bar G_{i})_{11}
\cdots (\bar G_{g})_{11}$, all elements $(M_i )_{11}$ and $(\bar G_{i})_{11}$
are invertible due to invertibility of $M_{11}$
 and one can easily
show that formulas (\ref{IS}) define the isomorphism of $\Lt$ and $\Lo
^{\otimes^g}$.

There is no natural anti-involution of the graph algebra. However, the
following anti-automorphism \cite{F} plays the role of $\ast$-operation on
$\Lg$:
\be 
\rho (A_i)=M_+A_i^{-1}M_+^{-1},\qquad \rho (B_i)=M_+B_i^{-1}M_+^{-1}
\label{aut}
\ee
The square of the anti-automorphism is equal to
\[ \rho^2 (A_i)=MA_iM^{-1},\qquad \rho^2 (B_i)=MB_iM^{-1}\]
It is worthwhile to note that this anti-automorphism acts on the elements
$M_{ij}$ as follows
\[ \rho (M_{\pm})=M_{\mp}\]

Let us introduce a set $\Phi$ of quantum constraints
\be
\Phi_{ij} =M_{ij} -\delta_{ij} 
\label{con}
\ee
where $M_{ij}$ are components of $M=M_1$.

Let us consider the left and right ideals of the graph algebra, generated by  
the set $\Phi$ and let ${\cal F}_{\I_{L}}$ (${\cal F}_{\I_{R}}$) be the maximal 
subalgebra 
of $\Lg$ such that
$\I_{L}$ ($\I_{R}$) is a two-sided ideal of ${\cal F}_{\I_{L}}$ 
(${\cal F}_{\I_{R}}$). Let $\F$ be the intersection of 
${\cal F}_{\I_{L}}$ and ${\cal F}_{\I_{R}}$:
\[ \F =\{ f\in \Lg : \I_{L} f\subset \I_{L}\quad and \quad f\I_{R}\subset 
\I_{R}\}.\]
It is obvious that $\F$ is a subalgebra of $\Lg$ and the elements 
$\Phi_{ij}\in \F$.

{\bf Definition 4.} Let $\I$ be the intersection of the algebra $\F$ with the
union of $\I_{L}$ and $\I_{R}$. The moduli algebra $\Mg$ is defined as a
quotient of the algebra $\F$ over the ideal $\I$
\be
\Mg = \F\big/ \I
\label{DM}
\ee
Some comments are in order. It is clear that $\I$ is a two-sided ideal of $\F$.
It follows from the action of the anti-automorphism $\rho$ on the elements
$M_{ij}$ that $\rho (\I_{L})=\I_{R}$ and $\rho (\I_{R})=\I_{L}$ and, hence if
$f\in\F$ and $i\in\I$ then $\rho (f)\in\F$ and $\rho (i)\in\I$. Thus the moduli
algebra inherits the anti-automorphism $\rho$ from $\Lg$. Moreover it is
possible to show that, being restricted to the moduli algebra, the
anti-automorphism becomes an anti-involution of $\Mg$. 

Our aim is to describe the structure of the centre ${\cal Z}(\Lg)$ of the graph 
algebra at roots of 1, and in the remainder of the paper we assume that $q$ is 
a primitive $p^{th}$ root of unity, p being odd. Due to the Proposition 1 we 
can begin with the study of the centre ${\cal Z}(\Lo)$ of the algebra $\Lo$.

\section{The centre of the graph algebra $\Lot$}

Let $a_{ij}$, $b_{ij}$ and $M_{11}^{-1}$ be the generators of $\Lot$ with the 
commutation
relations given by formulas (\ref{abcom}) from the Appendix.

{\bf Proposition 2.} 1) The centre ${\cal Z}(\Lo)$ of $\Lot$ is generated 
by the elements $a_{11}^{\pm p}$, $a_{12}^p$, $a_{21}^p$, $b_{11}^{\pm p}$, 
$b_{12}^p$,$b_{21}^p$ and $M_{11}^{-p}$, subject to the single relation
\be 
M_{11}^{-p}({\cal B}{\cal A}^{-1}{\cal B}^{-1}{\cal A})_{11}=1
\label{sr}
\ee
where ${\cal A}_{ij}=a_{ij}^p$, ${\cal B}_{ij}=b_{ij}^p$ 
for $i,j\neq 2$ simultaneously and $\det{\cal A}= \det{\cal B}=1$.

2) The algebra $\Lot$ is a free ${\cal Z}(\Lo)$-module with basis the set of 
monomials 
\[
a_{11}^{r_{1}}a_{12}^{s_{1}}a_{21}^{t_{1}}
b_{11}^{r_{2}}b_{12}^{s_{2}}b_{21}^{t_{2}}\quad {\mbox and}\quad 
0\leq r_{i},s_{i},t_{i}\leq p-1. \]

3) The ring $\Lot$ is an integral domain.

{\it Proof}. Let us introduce a new set of generators of $\Lot$ by means of the
following formulas
\[ X_1 =(a_{11}^{-1}a_{12}-b_{11}^{-1}b_{12})b_{11}^{2},\qquad 
X_2 =a_{11}a_{21}b_{11}^{-2}\]
\[ X_3 =a_{11}^{-2}b_{12}b_{11},\qquad 
X_4 =(b_{21}b_{11}^{-1}-a_{21}a_{11}^{-1})a_{11}^{2}\]
In terms of these generators the relations (\ref{abcom}) from the Appendix take 
the form
\bea
&&a_{11}b_{11}=qb_{11}a_{11},\qquad a_{11}X_i =X_i a_{11},
\qquad b_{11}X_i =X_i b_{11}\nonumber\\
&&X_1 X_4 =q^2 X_4 X_1,\qquad X_i X_{i+2}=q^{-2} X_{i+2}X_i ,\quad
i=1,2\nonumber\\
&&X_i X_{i+1} =q^2 X_{i+1}X_i +q^2 -1,\quad i=1,2,3\nonumber\\
&&M_{11}^{-1}X_1=q^{2}X_1M_{11}^{-1},\quad M_{11}^{-1}X_3=q^{2}X_3M_{11}^{-1}
\nonumber\\
&&M_{11}^{-1}X_2=q^{-2}X_2M_{11}^{-1},\quad M_{11}^{-1}X_4=q^{-2}X_4M_{11}^{-1}
\nonumber\\
&&M_{11}^{-1}q^{-2}(1+X_1 X_2 +X_1 X_4 +X_3 X_4 +X_1 X_2 X_3 X_4)=1
\label{newrel}
\eea
It is now of no problem to show that the elements $M_{11}^{-p}$, 
$a_{11}^{\pm p}$, 
$b_{11}^{\pm p}$ and $X_i^{p}$ 
generate the centre of $\Lot$. To do this one should use the following lemma,
which can be easily proved by induction

{\bf Lemma 1.} Let elements $c$, $Z$ and $W$ satisfy the relations $cZ=Zc$,
$cW=Wc$, $ZW=q^2 WZ +c$, then the following relation is valid
\[ Z^m W^n =\sum_{k=0}^{m} q^{2(m-k)(n-k)}c^k W^{n-k}Z^{m-k}\frac {(n)_q !}
{(n-k)_q !} \frac  {(m)_q !}{(m-k)_q !(k)_q !} \]
where $m\leq n$, $(n)_q =\frac {1-q^{2n}}{1-q^2}$.

The first part of the Proposition 2 follows now from the simple relations
between the elements $a_{ij}^p$, $b_{ij}^p$ and $X_{i}^p$
\[ X_1^p =(a_{11}^{-p}a_{12}^p -b_{11}^{-p}b_{12}^p )b_{11}^{2p},\qquad 
X_2^p =a_{11}^p a_{21}^p b_{11}^{-2p}\]
\[ X_3^p =a_{11}^{-2p}b_{12}^p b_{11}^p ,\qquad 
X_4^p =(b_{21}^p b_{11}^{-p}-a_{21}^p a_{11}^{-p})a_{11}^{2p}\]
which can be proved by using the well-known Lemma

{\bf Lemma 2.} Let elements $a$ and $b$ satisfy the relation $ab=q^2 ba$ or
$a(a+b)=q^2 (a+b)a$ and $q^p =1$, then
\[ (a+b)^p =a^p +b^p\]
The relation (\ref{sr}) will be proved later in this section (see Proposition
4).

The second part of the Proposition 2 follows from the obvious observation that
the products $a_{11}^{\pm r_{1}}a_{12}^{s_{1}}a_{21}^{t_{1}}
b_{11}^{\pm r_{2}}b_{12}^{s_{2}}b_{21}^{t_{2}}M_{11}^{-v}$, where  
$r_{i},s_{i},t_{i},v\in
N$, are a basis of $\Lot$.

Relations (\ref{newrel}) show that $\Lot$ is the tensor product of the Weil
algebra, generated by $a_{11}$ and $b_{11}$, and the algebra ${\cal X}$
generated by $X_i$ and $M_{11}^{-1}$ and, therefore, to prove that $\Lot$ is
an integral domain, it is enough to show that ${\cal X}$ is an integral
domain. We have to prove that if $fg=0$ then either $f=0$ or $g=0$.

\noindent An arbitrary element $f\in {\cal X}$ can be presented in the form
\[ f=\sum_{i,l=0}^{\infty}\sum_{j,k=0}^{p-1}\,
f_{ijkl}(M_{11}^{-p},X_2^p,X_3^p)X_1^iX_2^jX_3^kX_4^l\]
It is obvious that if $f\ne 0$, then $X_1f\ne 0$. Let us show that $fX_1= 0$
if and only if $f=0$. It is clear that it is enough to prove the statement 
only for elements of the form
\[ f=\sum\, f_{ij}X_1^iX_2^j\]
Using the commutation relation for $X_1$ and $X_2$ one gets
\bea
fX_1&=&\sum\, f_{ij}\left(
q^{-2j}X_1^{i+1}X_2^j+(q^{-2j}-1)X_1^{i}X_2^{j-1}\right)\nonumber\\
&=&\sum\, \left(
f_{ij}q^{-2j}+(q^{-2(j+1)}-1)f_{i+1,j+1}\right) X_1^{i+1}X_2^j
\nonumber
\eea
This expression is equal to zero only if
\[f_{ij}=(q^{2j}-q^{-2})f_{i+1,j+1}\]
Let $j=kp+j_0$, $0\le j_0\le p-1$. By simple induction one gets
\[ f_{ij}= (q^{2j_0}-q^{-2})(q^{2(j_0+1)}-q^{-2})\cdots 
(q^{2(j_0+p-j_0-1)}-q^{-2}) f_{i+p-j_0,(k+1)p}=0 \]
Thus if  $fX_1= 0$ then  $f=0$. In the same manner one can show that if 
$f\ne 0$ then 
$fX_4\ne 0$, $X_4f\ne 0$, $f(1+X_1X_2)\ne 0$, $(1+X_1X_2)f\ne 0$.

\noindent Let us note that for generic values of $q$ this statement is not
valid and the ring ${\cal X}$ is not an integral domain.

\noindent Let us consider elements 
$\tilde f =X_1^{p-1}X_4^{p-1}(1+X_1X_2)^{p-1}f$ and 
$\tilde g = gX_1^{p-1}X_4^{p-1}(1+X_1X_2)^{p-1}$. We have just shown that
equation $\tilde f =0$ ($\tilde g =0$) is equivalent to equation 
$f=0$ ($g=0$), therefore it is enough to show that if 
$\tilde f\tilde g =0$ then either $\tilde f =0$ or $\tilde g =0$. Elements 
$\tilde f$ and $\tilde g$ can be written in the form
\[\tilde f=\sum_{j,l=0}^{\infty}\sum_{i,k=0}^{p-1}\,
f_{ijkl}(M_{11}^{-p},X_2^p,X_3^p)M_{11}^iX_1^j(1+X_1X_2)^kX_4^l\]
\[\tilde g =\sum_{j,l=0}^{\infty}\sum_{i,k=0}^{p-1}\,
g_{ijkl}(M_{11}^{-p},X_2^p,X_3^p)M_{11}^iX_1^j(1+X_1X_2)^kX_4^l\]
The product of these elements is equal to
\bea
\tilde f\tilde g &=&\sum_{j_1,j_2,l_1,l_2=0}^{\infty}
\sum_{i_1,i_2,k_1,k_2=0}^{p-1}\,q^{2i_2(j_1-l_1)-2j_2(k_1+l_1)}f_{i_1j_1k_1l_1}
g_{i_2j_2k_2l_2}\times\nonumber\\
&& M_{11}^{i_1+i_2}X_1^{j_1+j_2}(1+X_1X_2)^{k_1+k_2}X_4^{l_1+l_2}
\nonumber
\eea
It is now of no problem to show that $\tilde f\tilde g =0$ 
if and only if either $\tilde f =0$ or $\tilde g =0$.

\noindent This proves the Proposition.

One can introduce the Poisson structure on the centre of $\Lot$ by means of the
following formula (see, for example, \cite{CP})
\be
\frac {k}{2\pi}\{ x,y\} = \lim _{\kappa \to q}\frac{xy-yx}{1-\kappa^{p^2}}
\label{PS}
\ee
Let $\PLg^*$ be the extension of $\PLg$ by means of the element $M_{11}^{-1}$. 
Using relations (\ref{newrel}) and Lemma 1 one can easily calculate the Poisson
brackets between the generators of ${\cal Z}(\Lo)$ and prove the following
proposition

{\bf Proposition 3.} The centre ${\cal Z}(\Lo)$ of $\Lot$ endowed with the
Poisson structure (\ref{PS}) is isomorphic to the extended Poisson graph algebra
 ${\cal PL}_1^*(sl_2)$
and the isomorphism is given by the formulas
\[ \phi (a_{ij}^p)=\a_{ij}, \qquad \phi (b_{ij}^p)=\b_{ij},\]
where $\a_{ij}$, $\b_{ij}$ are generators of the extended Poisson graph algebra 
${\cal PL}_1^*(sl_2)$
and $i,j\neq 2$ simultaneously.

{\bf Remark 2.} Let us note that as a by-product we have proven the well-known
theorem that ${\cal Z}_0 ({\cal U}_q (sl_2))$ is isomorphic to
$C[SL_2^* ]$ (see, e.g. \cite{DK,R,CP}).

To proceed with the study of the centre of the graph algebra $\Lg$ we will need
to know how the central elements  
${\cal M}_{ij}=M_{ij}^p =(BA^{-1}B^{-1}A)_{ij}^p$ are
expressed through the elements $a_{ij}^p$ and $b_{ij}^p$. Another reason to 
find these expressions is that
the ideal of ${\cal Z}(\Lg)$, generated by the elements $M_{ij}^p - 
\delta_{ij}$,
belongs to the centre ${\cal Z}({\cal I})$.

{\bf Proposition 4.} The central elements $M_{ij}^p$ are expressed through the
generators $a_{ij}^p$ and $b_{ij}^p$ of ${\cal Z}(\Lot)$ by means of the
following formula
\[ {\cal M}={\cal B}{\cal A}^{-1}{\cal B}^{-1}{\cal A}, \]
where ${\cal A}_{ij}=a_{ij}^p$, ${\cal B}_{ij}=b_{ij}^p$ 
for $i,j\neq 2$ simultaneously and $\det{\cal A}= \det{\cal B}=1$.

{\it Proof}. Let us introduce matrices $D=BA^{-1}$ and $C=B^{-1}A$. We firstly 
show that the
matrix ${\cal  M}$ is expressed through the matrix elements $d_{ij}^p$,
$c_{ij}^p$ as follows 
\[ {\cal  M}={\cal  D}{\cal  C}\]
where ${\cal  D}_{ij}=d_{ij}^p$, ${\cal  C}_{ij}=c_{ij}^p$ for $i,j\neq 2$
simultaneously.

The matrices $D$ and $C$ have the following commutation relations
\[ D^{1}\Rp D^{2}\Rpi =\Ra D^{2}\Rai D^{1},\quad
C^{1}\Rp C^{2}\Rpi =\Ra C^{2}\Rai C^{1}\]
\[ C^{1}\Rp D^{2}\Rpi =\Rp D^{2}\Rpi C^{1},\qquad 
\detq C =\detq C =q^3 \]
which are rewritten in components in the Appendix.

\noindent Using the relations (\ref{cdcom}) and Lemma 2 one gets
\[ M_{11}^p =(DC)^p_{11} =d_{11}^p(c_{11}+d_{11}^{-1}d_{12}c_{21})^p =
d_{11}^p c_{11}^p+d_{12}^p c_{21}^p =({\cal  D}{\cal  C})_{11}\]
\bea 
M_{12}^p &=&(DC)^p_{12} =(d_{11}c_{12}+q^2 d_{12}c_{11}^{-1}c_{21}c_{12}+
q^3 d_{12}c_{11}^{-1})^p \nonumber\\
&=&(d_{11}c_{12}+q^2 d_{12}c_{11}^{-1}c_{21}c_{12})^p +(d_{12}c_{11}^{-1})^p =
(d_{11}+q^2 d_{12}c_{11}^{-1}c_{21})^p c_{12}^p +d_{12}^p c_{11}^{-p} 
\nonumber\\
&=&d_{11}^p c_{12}^p +d_{12}^p c_{11}^{-p}c_{21}^p c_{12}^p +
d_{12}^p c_{11}^{-p} =({\cal  D}{\cal  C})_{12}
\nonumber
\eea
\bea
M_{21}^p &=&(DC)^p_{21} =(d_{21}c_{11}+q^2 d_{11}^{-1}d_{21}d_{12}c_{21}+
q^3 d_{11}^{-1} c_{21})^p \nonumber\\
&=&d_{11}^{-p} c_{21}^p + (d_{21}c_{11}+q^2 d_{11}^{-1}d_{21}d_{12}c_{21})^p =
d_{11}^{-p} c_{21}^p + d_{21}^p (c_{11}+ d_{11}^{-1}d_{12}c_{21})^p 
\nonumber\\
&=&d_{11}^{-p} c_{21}^p + d_{21}^p (c_{11}^p + d_{11}^{-p}d_{12}^p c_{21}^p) 
 =({\cal  D}{\cal  C})_{21}
\nonumber
\eea
To complete the proof of the Proposition 4 one should show that 
\[
{\cal  C}={\cal  B}^{-1}{\cal  A},\quad {\cal  D}={\cal  B}{\cal  A}^{-1}
\]
It can be done by using the following lemma which can be easily proved
with the help of the Lemma 2:

{\bf Lemma 3.} Let elements $c$, $Z$ and $W$ satisfy the relations $cZ=Zc$,
$cW=Wc$, $ZW-c=q^{-2}( WZ -c)$ and $q^p =1$ then 
\[ (ZW-c)^p =Z^p W^p -c^p \]
Then the calculation of $c_{ij}^p$ gives
\bea
c_{11}^p &=&(B^{-1}A)_{11}^p =\left( (q^2 b_{22} +(1-q^2 )b_{11})a_{11} - 
q^2 b_{12}a_{21}\right)^p \nonumber\\
&=&\left( q^2 b_{11}^{-1}a_{11} +q^4 b_{11}^{-1}b_{21}b_{12}a_{11} +
 (1-q^2 )b_{11}a_{11} -q^2 b_{12}a_{21}\right)^p \nonumber\\
&=&\left( b_{11}^{-1}a_{11} +q^2 b_{11}^{-1}b_{21}b_{12}a_{11}
-qa_{21}b_{12}\right)^p 
\nonumber\\
&=&b_{11}^{-p}a_{11}^p \left( 1 +q^2 (b_{21} -a_{21}a_{11}^{-1}b_{11})b_{12} 
\right)^p 
\nonumber\\
&=&b_{11}^{-p}a_{11}^p \left( 1 + (b_{21} -a_{21}a_{11}^{-1}b_{11})^p b_{12}^p
\right) =
({\cal  B}^{-1}{\cal  A})_{11}
\nonumber
\eea
\bea
c_{12}^p &=&(B^{-1}A)_{12}^p =(q^2 b_{22}a_{12} - qa_{22}b_{12})^p 
\nonumber\\
&=&(-a_{11}^{-1}b_{12} -q^2 a_{11}^{-1}a_{21}a_{12}b_{12} +q b_{22}a_{12})^p
 \nonumber\\
&=&-a_{11}^{-p}b_{12}^p + (-q^2 a_{11}^{-1}a_{21}a_{12}b_{12} +q b_{22}a_{12})^p
 \nonumber\\
&=&-a_{11}^{-p}b_{12}^p +a_{12}^p (-qa_{11}^{-1}a_{21}b_{12} + b_{22})^p
 \nonumber\\
&=&-a_{11}^{-p}b_{12}^p +a_{12}^p \left( b_{11}^{-1}+(b_{21}b_{11}^{-1}-
a_{21}a_{11}^{-1})b_{12}\right)^p
 \nonumber\\ 
&=&-a_{11}^{-p}b_{12}^p +a_{12}^p b_{11}^{-p}\left( 1+(b_{21}b_{11}^{-1}-
a_{21}a_{11}^{-1})b_{12}b_{11}\right)^p
 \nonumber\\ 
&=&-a_{11}^{-p}b_{12}^p +a_{12}^p b_{11}^{-p}\left( 1+(b_{21}-
a_{21}a_{11}^{-1}b_{11})^p b_{12}^p\right) =({\cal  B}^{-1}{\cal  A})_{12}
\nonumber
\eea
\bea
c_{21}^p &=&(B^{-1}A)_{21}^p =(-q^2 b_{21}a_{11} + b_{11}a_{21})^p \nonumber\\
&=&a_{11}^p (-q^2 b_{21} + b_{11}a_{21}a_{11}^{-1})^p =
a_{11}^p (- b_{21}^p + b_{11}^p a_{21}^p a_{11}^{-p}) =
({\cal  B}^{-1}{\cal  A})_{21}
\nonumber
\eea
The calculation of the matrix elements ${\cal D}_{ij}$ can be done in
the same manner.

The matrix elements of $M$ can be expressed through the generators $a_{11}$,
$b_{11}$ and $X_i$ as follows
\bea
M_{11}&=&1+X_2 X_1 +X_1 X_4 +X_3 X_4 +X_3 X_2 X_1 X_4\nonumber\\
&=&q^{-2}(1+X_1 X_2 +X_1 X_4 +X_3 X_4 +X_1 X_2 X_3 X_4)\nonumber\\
M_{12}&=&-q^{-2}((1+X_1 X_2)X_3 +X_1) +q^{-2}X_1 M_{11}b_{11}^{-2}+
X_3 M_{11}a_{11}^2 b_{11}^{-2}\nonumber\\
M_{21}&=&-q^{-2}(X_2 (1+X_3 X_4) +X_4) +X_4 M_{11}a_{11}^{-2}+
q^2 X_2 M_{11}a_{11}^{-2} b_{11}^{2}
\label{20}
\eea
It follows from the Proposition 4 that
\bea
M_{11}^p &=&1+X_1^p X_2^p +X_1^p X_4^p +X_3^p X_4^p +X_1^p X_2^p X_3^p X_4^p
\nonumber\\
M_{12}^p &=&-(1+X_1^p X_2^p )X_3^p -X_1^p +X_1^p M_{11}^p b_{11}^{-2p}+
X_3^p M_{11}^p a_{11}^{2p} b_{11}^{-2p}\nonumber\\
M_{21}^p &=&-X_2^p (1+X_3^p X_4^p) -X_4^p +X_4^p M_{11}^p a_{11}^{-2p}+
X_2^p M_{11}^p a_{11}^{-2p} b_{11}^{2p}
\label{21}
\eea
One can use eqs.(\ref{20}) to prove the Proposition 4.

We are now ready to discuss the structure of the centre of the graph algebra
$\Lg$ and the moduli algebra $\Mg$.

\section{The centre of the graph and moduli algebras}

{\bf Proposition 5.} 1) The centre ${\cal Z}(\Lg )$ of $\Lt$ is generated
by the elements $a_{i11}^{\pm p}$, $a_{i12}^p$, $a_{i21}^p$, 
$b_{i11}^{\pm p}$, $b_{i12}^p$,$b_{i21}^p$ and $M_{11}^{-p}$, subject to the
single relation
\[ M_{11}^{-p}({\cal B}_g{\cal A}_g^{-1}{\cal B}_g^{-1}{\cal A}_g\cdots 
{\cal B}_1{\cal A}_1^{-1}{\cal B}_1^{-1}{\cal A}_1)_{11}=1 \]

2) The algebra $\Lt$ is a free ${\cal Z}(\Lot)$-module with basis the set of 
monomials 
$\prod_{i=1}^g \left( a_{i11}^{r_{i1}}a_{i12}^{s_{i1}}a_{i21}^{t_{i1}}
b_{i11}^{r_{i2}}b_{i12}^{s_{i2}}b_{i21}^{t_{i2}}\right)$ and 
$0\leq r,s,t\leq p-1$.

3) The ring $\Lt$ is an integral domain.

4) The centre  ${\cal Z}(\Lg )$ endowed with the Poisson structure (\ref{PS})
is isomorphic to the extended Poisson graph algebra ${\cal PL}_g^*(sl_2)$.

5) The centre ${\cal Z}(\Lg )$ is isomorphic to the tensor product of $g$
copies of the extended Poisson graph algebra ${\cal PL}_g^*(sl_2)$. The 
isomorphism is given by means of the following formulas
\be
{\cal A}_i ={\cal M}_{+}^{-1}(i+1)\bar {\cal A}_i {\cal M}_{+}(i+1),\qquad 
{\cal B}_i ={\cal M}_{+}^{-1}(i+1)\bar {\cal B}_i {\cal M}_{+}(i+1)
\label{is}
\ee
where
\bea
({\cal A}_i)_{mn} &=&(A_i )_{mn}^p ,\quad ({\cal B}_i)_{mn} =(B_i )_{mn}^p
\nonumber\\ 
(\bar {\cal A}_i)_{mn} &=&(\bar A_i )_{mn}^p ,\quad (\bar {\cal B}_i)_{mn} 
=(\bar B_i )_{mn}^p \nonumber
\eea
and $m, n\neq 2$ simultaneously,
\be 
{\cal M}_{i}={\cal M}_{-}^{-1}(i){\cal M}_{+}(i)=
{\cal B}_g {\cal A}_g^{-1}{\cal B}_g^{-1} {\cal A}_g\cdots 
{\cal B}_i {\cal A}_i^{-1}{\cal B}_i^{-1} {\cal A}_i
\label{Mi}
\ee
\bea
{\cal M}_{\pm}(i)&=&\bar {\cal G}_{\pm}(i)\bar {\cal G}_{\pm}(i+1)\cdots
\bar {\cal G}_{\pm}(g)\nonumber\\ 
\bar{\cal  G}_{i}&=&\bar {\cal G}_-^{-1}(i)\bar {\cal G}_+(i)
=\bar {\cal B}_i\bar {\cal A}_i^{-1}\bar {\cal B}_i^{-1}
\bar {\cal A}_i
\label{Mpm}
\eea
and $\bar A_i$, $\bar B_i$ are matrices from the Proposition 1.

{\it Proof}. It is obvious that Propositions 5.1-5.4  follow from
Propositions 1, 2, 3, and 5.5, thus it is enough to prove the Proposition 5.5.
The matrix elements of $M_+ (i)$ can be expressed through the matrix elements
of $M_i$ as follows
\bea
M_+ (i)_{11} &=& m_{i11}^{\frac{1}{2}},\quad 
M_+ (i)_{12} = m_{i11}^{-\frac{1}{2}}m_{i12},\quad 
M_+ (i)_{22} =m_{i11}^{-\frac{1}{2}},\nonumber\\  
M_+^{-1} (i)_{11} &=& m_{i11}^{-\frac{1}{2}},\quad 
M_+^{-1} (i)_{12} =-q m_{i11}^{-\frac{1}{2}}m_{i12},\quad 
M_+^{-1} (i)_{22} =m_{i11}^{\frac{1}{2}}
\nonumber
\eea
Using formulas (\ref{IS}) one gets
\[ 
(A_i )_{11}^p =(\bar a_{i11} -q^2 m_{i+1,12}\bar a_{i21})^p = \bar a_{i11}^p -
m_{i+1,12}^p \bar a_{i21}^p =
(\tilde {\cal M}_{+}^{-1}(i+1)\bar {\cal A}_i \tilde {\cal M}_{+}(i+1))_{11}\]
\[ (A_i )_{21}^p =(m_{i+1,11}\bar a_{i21})^p = m_{i+1,11}^p \bar a_{i21}^p =
(\tilde {\cal M}_{+}^{-1}(i+1)\bar {\cal A}_i \tilde {\cal M}_{+}(i+1))_{21}\]
\bea
&&(A_i )_{12}^p 
=(\bar a_{i11}m_{i+1,11}^{-1}m_{i+1,12} + 
\bar a_{i12}m_{i+1,11}^{-1} - \bar a_{i21}m_{i+1,11}^{-1}m_{i+1,12}^2 - 
\bar a_{i22}m_{i+1,11}^{-1}m_{i+1,12})^p \nonumber\\
&& =\bar a_{i11}^{-p} m_{i+1,11}^{-p}(\bar a_{i11}^2 m_{i+1,12} + 
\bar a_{i11}\bar a_{i12} - \bar a_{i11}\bar a_{i21}m_{i+1,12}^2  
- m_{i+1,12} - q^2 \bar a_{i21}\bar a_{i12}m_{i+1,12})^p \nonumber\\
&&=\bar a_{i11}^{-p} m_{i+1,11}^{-p}\left( 
(\bar a_{i11} -q^2 \bar a_{i21}m_{i+1,12})
(\bar a_{i11}m_{i+1,12}+ \bar a_{i12})-m_{i+1,12}\right) ^p \nonumber\\
&&=\bar a_{i11}^{-p} m_{i+1,11}^{-p}\left( 
(\bar a_{i11} -q^2 \bar a_{i21}m_{i+1,12})^p
(\bar a_{i11}m_{i+1,12}+ \bar a_{i12})^p -m_{i+1,12}^p\right) \nonumber\\
&&=(\tilde {\cal M}_{+}^{-1}(i+1)\bar {\cal A}_i \tilde {\cal M}_{+}(i+1))_{12}
\label{3}
\eea
where $\tilde {\cal M}_{i+1,kl}=m_{i+1,kl}^p$.

To get these expressions one should use the commutativity of $m_{i+1,kl}$ and 
$ \bar a_{imn}$, the Lemmas 2 and 3 and 
eq.(\ref{3}) follows from the identification
\[ \bar a_{i11} -q^2 \bar a_{i21}m_{i+1,12} =Z,\quad 
\bar a_{i11}m_{i+1,12}+ \bar a_{i12} =W,\quad m_{i+1,12} = c \]

One sees that to prove the Proposition 5.4 one should show that 
\be
\tilde {\cal M}_{i} ={\cal M}_{i}
\label{4}
\ee
The matrix elements of $M_i$ can be expressed through the matrix elements of 
$\bar G_i$ as follows
\[ m_{i11}= \bar g_{i11}\bar g_{i+1,11}\cdots\bar g_{g11}\]
\[ m_{i12}= \sum_{k=i}^{g}\,\bar g_{i11}\cdots\bar g_{k-1,11}\bar g_{k12}\]
\[ m_{i21}= \sum_{k=i}^{g}\,\bar g_{i11}\cdots\bar g_{k-1,11}\bar g_{k21}\]

Using these expressions and the Lemma 2 one immediately gets that 
\[ \tilde {\cal M}_{+}(i) = \tilde\Gamma_{\pm} (i)\cdots \tilde\Gamma_{\pm} (g)
\]
where $(\tilde\Gamma_{\pm} (i))_{mn} = (\bar {\cal G}_{\pm}(i))_{mn}^p$.

Thus to complete the proof of the Proposition 5.4 it remains to remember that
the equation
\[ \tilde\Gamma_{i} = \tilde\Gamma_{-}^{-1}(i)\tilde\Gamma_{+} (i) =
\bar{\cal  G}_{i}=\bar {\cal B}_i\bar {\cal A}_i^{-1}\bar {\cal B}_i^{-1}
\bar {\cal A}_i\]
was proved in Proposition 4.

\noindent The Proposition 5 is proved.

{\bf Remark} 3. Strictly speaking to define the matrices ${\cal M}_{\pm}(i)$, 
$\bar {\cal G}_{\pm}(i)$ one should use the elements $(m_{i11}^{\pm
\frac{1}{2}})^p$, $(\bar g_{i11}^{\pm \frac{1}{2}})^p$ from the extended graph
algebras. However the formulas (\ref{is}) (and (\ref{IS})) do not depend on
them due to the obvious relations $(m_{i11}^{\frac{1}{2}})^2 =m_{i11}$, 
$(\bar g_{i11}^{\frac{1}{2}})^2 =\bar g_{i11}$ and we use the matrices 
${\cal M}_{\pm}$, $\bar {\cal G}_{\pm}$ only to simplify the notations and the
proof of eq.(\ref{4}).

{\bf Remark 4.} The method described in the third section to prove the
Proposition 4 can be applied to prove eq.(\ref{Mi}) 
without using the decomposion (\ref{Mpm}).

Let us proceed with the study of the centre ${\cal Z}(\Mg)$ of the moduli
algebra. 

{\bf Proposition 6.} The centre ${\cal Z}(\F)$ coincides with ${\cal Z}(\Lg)$  
and
the centre ${\cal Z}(\I)$ is the ideal of  ${\cal Z}(\Lg)$, generated by the
elements $M_{ij}^p -\delta_{ij}$.

{\it Proof}. It is obvious that ${\cal Z}(\Lg)\subset {\cal Z}(\F)$. As
was proved in Proposition 5, $\Lg$ and, hence $\F$ are integral domains and,
therefore ${\cal Z}(\I)\subset {\cal Z}(\F)$. 
Let $z$ belongs to ${\cal Z}(\F)$. It is not difficult to show that the
elements $\bar a_{k11}M_{21}^{p-1}M_{12}^{p-1}$ and 
$\bar b_{k11}M_{21}^{p-1}M_{12}^{p-1}$ belong to $\F$. The ring $\Lg$ is an
integral domain and, hence the following equations should be valid
\[ z\bar a_{k11}=\bar a_{k11}z, \qquad z\bar b_{k11}=\bar b_{k11}z \]
It follows from these equations that the element $z$ should be of the form 
\[ z=z((X_{k})_i, \bar a_{k11}^p, \bar b_{k11}^p)\].
Taking into account that $z$
commutes with $M_{12}$ and $M_{21}$ one gets
\[ z(X_{k})_i=(X_{k})_iz \quad \forall k=1,...,g\quad and\quad \forall
i=1,2,3,4.\]
Therefore, $z$ belongs to ${\cal Z}(\Lg)$ and ${\cal Z}(\F)$ 
coincides with ${\cal Z}(\Lg)$.

It can be easily shown using eqs.(\ref{20}) and (\ref{21}) that 
any element
$z\in {\cal Z}(\I)$ belongs to the ideal of ${\cal Z}(\Lg)$ generated by the
elements  $M_{ij}^p -\delta_{ij}$.

Let us consider a subalgebra ${\cal J}$ of $\F$ which consists of the elements
such that the commutator of an element $j\in {\cal J}$ with any element
$f\in\F$ belongs to the ideal $\I$:
\[ {\cal J}=\{ j\in \F : jf-fj\in \I \quad\forall f\in\F\}\]
It is obvious that ${\cal Z}(\Lg)\subset {\cal J}$ is the centre of ${\cal J}$,
 $\I$ is an ideal of 
${\cal J}$ and the centre ${\cal Z}(\Mg)$ coincides with the factor algebra of 
${\cal J}$ over $\I$
\[ {\cal Z}(\Mg)= {\cal J}/\I \]
As was shown before the graph algebra $\Lg$ and, therefore, the algebras $\F$
and ${\cal J}$ are finitely generated over ${\cal Z}(\Lg)$. As has been just
proved ${\cal Z}(\I)\subset {\cal Z}(\Lg)$ and, hence the moduli algebra 
$\Mg =\F/\I$ and the centre ${\cal Z}(\Mg)$ are finitely generated over 
${\cal Z}(\Lg)/{\cal Z}(\I)$. Let us consider $Spec({\cal Z}(\Mg))$ and
$Spec({\cal Z}(\Lg)/{\cal Z}(\I))$, i.e. the set of all algebra homomorphisms
from ${\cal Z}(\Mg)$ and ${\cal Z}(\Lg)/{\cal Z}(\I)$ to $C$. It is clear that 
$Spec({\cal Z}(\Lg))$ is isomorphic to $C^{4g}\times (C^{\times})^{2g}$, i.e. a
complex affine space of dimension $6g$ with $2g$ hyperplanes of codimension 1
removed. Then $Spec({\cal Z}(\Lg)/{\cal Z}(\I))$ is a submanifold of 
$C^{4g}\times (C^{\times})^{2g}$ which is singled out by means of eq.(\ref{R}).

\section{Irreducible representations of the graph algebra}

It is clear that every irreducible $\Lg$-module $V$ is finite dimensional and
the centre ${\cal Z}(\Lg)$ acts by scalar operators on $V$ and, therefore,
there is a homomorphism $\chi_V \in Spec({\cal Z}(\Lg)): {\cal
Z}(\Lg)\rightarrow C$, the central character of $V$, such that
\[ z.v=\chi_V (z)v \]
for all $z\in {\cal Z}(\Lg)$ and $v\in V$.

Let us consider an ideal $\I_{\chi}$ of $\Lg$ generated by elements 
$z-\chi_V (z)$, $z\in {\cal Z}(\Lg)$. Then every irreducible representation of
the algebra $\Lg^{\chi} =\Lg/\I_{\chi}$ is an irreducible representation of
$\Lg$ with the central character $\chi_V (z)$.

Due to the Proposition 1 it is enough to construct representations of $\Lo$.
As was shown in the third section the algebra $\Lo$ is isomorphic to the tensor
product of the Weil algebra, generated by $a_{11}$ and $b_{11}$, and the
algebra ${\cal X}$ generated by $X_i$. There is no problem in constructing
representations of the Weil algebra and we begin with the discussion of
irreducible representations of ${\cal X}$.

{\bf Proposition 7.} The algebra  ${\cal X}_{\chi}={\cal X}/\I_{\chi}$, $\chi
(M_{11}^p)\ne 0$ is a simple $p^4$-dimensional algebra and, hence, the algebra
 $\Lo^{\chi}$ is simple $p^6$-dimensional.
 
{\it Proof}. We are going to show that the unity element belongs to an
arbitrary (nonzero) ideal ${\cal J}$ and, hence, the ideal coincides with 
${\cal X}_{\chi}$ and ${\cal X}_{\chi}$ is a simple algebra.

An arbitrary element $f\in {\cal X}_{\chi}$ can be presented in the form
\bea
f&=&\sum \,c_{ijkl}X_1^i (1+X_1X_2)^j(1+X_3X_4)^kX_4^l + 
 d_{ijkl}X_1^i (1+X_1X_2)^j(1+X_3X_4)^kX_3^l \nonumber\\
 &+&f_{ijkl}X_2^i (1+X_1X_2)^j(1+X_3X_4)^kX_4^l + 
 g_{ijkl}X_2^i (1+X_1X_2)^j(1+X_3X_4)^kX_3^l \nonumber
 \eea
where $i,j,k,l=0,1,...,p-1$.

\noindent Let $f$ belong to an ideal ${\cal J}$ and  let $l_2 (l_3)$ be the
maximal power of $X_2 (X_3)$ in the element $f$. Then one can easily show that
the element $X_1^{l_2}fX_4^{l_3}$, which obviously belongs to ${\cal J}$, can
be represented in the form
\[X_1^{l_2}fX_4^{l_3}= \sum \,c_{ijkl}X_1^i (1+X_1X_2)^j(1+X_3X_4)^kX_4^l \]
with some new coefficients $c_{ijkl}$.

\noindent This element can be rewritten as follows
\[X_1^{l_2}fX_4^{l_3}= \sum \,c_{\alpha ijk}M_{11}^{\alpha} (1+X_1X_2)^i
X_1^jX_4^k + d_{\alpha ijk}M_{11}^{\alpha} (1+X_3X_4)^i X_1^jX_4^k \]
where $M_{11}$ is given by eq.(\ref{20}).

\noindent Multiplying  $X_1^{l_2}fX_4^{l_3}$ on $(1+X_1X_2)^l$, where $l$ is
the maximal power of $(1+X_3X_4)$ in  $X_1^{l_2}fX_4^{l_3}$ one gets that the
ideal  ${\cal J}$ contains an element of the form
\[ \sum \,c_{\alpha ijk}M_{11}^{\alpha} (1+X_1X_2)^iX_1^jX_4^k\]

Let us now suppose that $X_1^pX_4^p\ne 0$. Then the elements $X_1$ and $X_4$
are invertible and using the commutation relations of $M_{11}$, $X_1$ and $X_4$
with $X_i$ one can easily show that the element of the form
\[ (1+X_1X_2)^l\sum \,c_{\alpha} M_{11}^{\alpha}X_1^{p-\alpha} X_4^{p-\alpha}
\]
belongs to ${\cal J}$.

\noindent If $X_2^p = 0$ then the element  $(1+X_1X_2)$ is invertible and from
the commutation relations of $(1+X_1X_2)$ with  $X_i$ one gets that the ideal 
${\cal J}$ contains $M_{11}$ and, hence the unity element.

\noindent If $X_2^p\ne 0$ then using the commutation relations of $X_2$ with
$X_i$ one gets that the element  $(1+X_1X_2)^l$ belongs to  ${\cal J}$. Now
using the relation
\[ (1+X_1X_2)(1+X_3X_4)=(1+X_3X_4)(1+X_1X_2) +(q^2 -1)X_1X_4 \]
one can easily show that  ${\cal J}$ contains the element $X_1X_4$ and,
therefore, the unity element.

Let us now consider the case $X_1^p =0$. Then the element  $(1+X_1X_2)$ is 
invertible and from the commutation relations of $M_{11}$ and $(1+X_1X_2)$ with
$X_i$ one can get that the element of the form
\[ X_1^{l_1}X_4^{l_4}\sum \,c_{\alpha i}M_{11}^{\alpha} (1+X_1X_2)^i \]
belongs to ${\cal J}$.

Let  $X_4^p\ne 0$. Using the commutation relations of $X_4$ with
$X_i$ one gets that ${\cal J}$ contains the element of the form
\[ X_1^{l_1}\sum \,c_{i}(1+X_1X_2)^i \]
With the help of $X_2$ one gets that the element
$ X_1^{l_1}X_2^{l_2}(1+X_1X_2)^l $ and, hence, $ X_1^{l_1}X_2^{l_2}$ belongs to
${\cal J}$ and using firstly $X_1$ and then $X_2$ one sees that the unity
element is in ${\cal J}$.

If $X_4^p =0$ then the element $(1+X_3X_4)$ is invertible and using
this element and $X_2$ one shows that the element 
$X_4^{p-1}\sum \,c_{\alpha }M_{11}^{\alpha} $ belongs to ${\cal J}$.
Using again $X_2$ one gets that ${\cal J}$ contains the element 
$X_4^{p-1}X_2^l$. With the help of $X_1$ and then $X_3$ one gets that the unity
element belongs to  ${\cal J}$.
The case $X_4^p =0$ , $X_1\ne 0$ can be considered in the same manner.

\noindent The Proposition 7 is proved.

The algebra ${\cal X}_{\chi}$, being simple, is isomorphic to $M_{p^2}(C)$ and,
therefore, the following proposition is valid.

{\bf Proposition 7.} Every irreducible representation of  ${\cal X}$ is
isomorphic to one of the following $p^2$-dimensional:
\bea
1)&& X_1\Psi (k,l)=\Psi (k+1,l)\nonumber\\
&& X_2\Psi (k,l)=q^{-2k}(y_2x_1+1-q^{2k})\Psi (k-1,l) + 
z_2 q^{-2(k+l)}\Psi (k,l+1)\nonumber\\
&& X_3\Psi (k,l)=q^{2(k+l)}(y_3x_4+1-q^{-2l})\Psi (k,l-1) + 
z_3 q^{2k}\Psi (k+1,l)\nonumber\\
&& X_4\Psi (k,l)=q^{-2k}\Psi (k,l+1)
\label{rep1}
\eea
Here $\Psi (k,l)$ is a basis of a $p^2$-dimensional vector space, 
$k,l=0,1,...,p-1$ and we use the following notations 
\[ \Psi (p,l)=x_1\Psi (0,l),\qquad x_1\Psi (-1,l)=\Psi (p-1,l)\]
\[ \Psi (k,p)=x_4\Psi (k,0),\qquad x_4\Psi (k,-1)=\Psi (k,p-1)\]
The complex parameters $x_1$, $x_4$, $y_2$, $y_3$,  $z_2$, $z_3$ should satisfy
the following equations
\[ (y_2x_1+1)(y_3x_4+1)\ne 0,\qquad z_2z_3=0 \]
\[ 1 +q^{-2}z_3(y_2x_1+1)+q^2z_2(y_3x_4+1)=0 \]
The central character of the representation is defined by the formulas
\[ \chi (X_1^p)=x_1,\qquad \chi (X_4^p)=x_4 \]
\[ \chi (X_2^p)=\frac {1}{x_1} ((y_2x_1+1)^p-1)+z_2^px_4 \]
\[ \chi (X_3^p)=\frac {1}{x_4} ((y_3x_4+1)^p-1)+z_3^px_1 \]
The element $M_{11}$ acts in the representation as follows
\[ M_{11}\Psi (k,l)=q^{2(l-k)}(y_2x_1+1)(y_3x_4+1)\Psi (k,l) \]
\bea
2)&& X_1\Psi (k,l)=\Psi (k+1,l)\nonumber\\
&& X_2\Psi (k,l)=-\Psi (k-1,l) + b_1q^{-2(k+l)}\Psi (k,l+1)
+ b_2q^{-2(k+2l)}\Psi (k+1,l+2)\nonumber\\
&& X_3\Psi (k,l)=-q^{2k}\Psi (k,l-1) + 
c_{p-1} q^{2(k+2l)}\Psi (p+k-1,p+l-2)\nonumber\\
&& X_4\Psi (k,l)=q^{-2k}\Psi (k,l+1)
\nonumber
\eea
where
\[ x_1x_4b_1b_2c_{p-1}\ne 0,\qquad 1+q^{10}x_1x_4b_2c_{p-1}=0 \]
\[ \chi (X_1^p)=x_1,\qquad \chi (X_4^p)=x_4 \]
\[ \chi (X_2^p)=(b_1^p+b_2^px_1x_4)x_4-x_1^{-1} \]
\[ \chi (X_3^p)=c_{p-1}^px_1^{p-1}x_4^{p-2}-x_4^{-1} \]
\[  M_{11}\Psi (k,l)=q^{2(l-k +3)}x_1x_4b_1c_{p-1}\Psi (k,l) \]

Now irreducible representations of the graph algebra $\Lot$ can be constructed
using, for example, the following representation of the Weil algebra
\[ a_{11}\Psi (m)=\Psi (m+1),\qquad b_{11}\Psi (m)=\beta_{11} q^{-m}\Psi (m)\]
\[\Psi (m+p)=\alpha_{11}\Psi (m),\quad \chi (a_{11}^p)=\alpha_{11},\quad
\chi (b_{11}^p)=\beta_{11}^p\]
Then the graph algebra acts in the tensor product of the representations of the
Weil algebra and the algebra ${\cal X}$
\bea
&&a_{11}\Psi (k,l,m)=\Psi (k,l,m+1)\nonumber\\
&&b_{11}\Psi (k,l,m)=\beta_{11} q^{-m}\Psi (k,l,m)\nonumber\\
&&a_{12}\Psi (k,l,m)=\beta_{11}^{-2} q^{2m}\Psi (k+1,l,m+1)+
\beta_{11}^{-2} z_3 q^{2(m+k+1)}\Psi (k+1,l,m+3)\nonumber\\
&&\qquad\qquad\qquad +\beta_{11}^{-2} q^{2(m+k+l+1)}(y_3x_4+1-q^{-2l})\Psi 
(k,l-1,m+3)
\nonumber\\
&&a_{21}\Psi (k,l,m)=\beta_{11}^{2} q^{2(m-k)}(y_2x_1+1-q^{2k})\Psi (k-1,l,m-1)
\nonumber\\&&\qquad\qquad\qquad +
\beta_{11}^{2} z_2 q^{2(m-k-l)}\Psi (k,l+1,m-1)\nonumber\\
&&b_{12}\Psi (k,l,m)=\beta_{11}^{-1} q^{m+2(k+l)}(y_3x_4+1-q^{-2l})
\Psi (k,l-1,m+2)\nonumber\\&&\qquad\qquad\qquad +
\beta_{11}^{-1} z_3 q^{m+2k}\Psi (k+1,l,m+2)\nonumber\\
&&b_{21}\Psi (k,l,m)=\beta_{11}^{3} q^{3m-2k+2}(y_2x_1+1-q^{2k})\Psi (k-1,l,m-2)
\nonumber\\&&\qquad\qquad\qquad +
\beta_{11}^{3} q^{3m-2(k+l)+2}\Psi (k,l+1,m-2) \nonumber\\&&\qquad\qquad\qquad +
\beta_{11} q^{m-2k}\Psi (k,l+1,m-2)
\nonumber
\eea
and we present formulas only for the representation (\ref{rep1}) of ${\cal X}$.

Using this representation and Proposition 1 one can easily construct all
irreducible representations of the graph algebra $\Lg$.

Let us briefly discuss representations of the moduli algebra $\Mg$. In this
case one should consider only representations with $\chi (M_{11}^p)=1$, 
$\chi (M_{12}^p)=\chi (M_{21}^p)=0$. Let $V_0$ be the submodule of an
irreducible left $\Lg$-module $V$ which is annihilated by the elements 
$ \Phi_{ij} =M_{ij} -\delta_{ij}$:
\[ V_0=\{ \Psi\in V:\quad \Phi_{ij}\Psi =0\} \]
and let $\delta V_0$ be a submodule of $V_0$ which consists of the vectors of
the following form:
\[ \delta V_0 = \{ \Psi\in V_0:\quad \Psi=\Phi_{ij}\chi_{ij}\quad
for\quad some\quad\chi_{ij}\in V \} \]
Then it is obvious that the moduli algebra acts in the factor module
$V_{ph}=V_0\big/\delta V_0 $.

\section{Conclusion}

In this paper we studied the structure of the centre and irreducible
representations of the graph algebra. The next problem to be solved is to
construct unitary representations of the graph algebra. The anti-automorphism
$\rho$ (\ref{aut}) can be used to define unitary representations of $\Lg$.
Namely, let $V$ be a left $\Lg$-module and $\beta$ be a bilinear form on
$V\times V$, such that $\beta (v_2,\l v_1)=\l \beta (v_2,v_1)$, 
$\beta (\l v_2, v_1)=\l^* \beta (v_2,v_1)$, $\l\in C$. A representation is
called unitary if $\beta (v_2,fv_1)=\beta (\rho (f)v_2,v_1)$ for all
$v_1,v_2\in V$ and $f\in \Lg$. It is obvious  that only representations with
central characters, satisfying the equations 
${\cal A}_i{\cal M}={\cal M}{\cal A}_i$, ${\cal B}_i{\cal M}={\cal M}{\cal
B}_i$, can be unitary.

We didn't study irreducible representations of the moduli algebra, however
there are some indications that $V_{ph}$ is an irreducible $\Mg$-module. One
should prove (or disprove) this conjecture and show that the dimension of
$V_{ph}$ is given by Verlinde's formula.

It would be interesting to clarify the relation between this approach to
quantization of the moduli space and the geometric quantization \cite{ADW,JW}. 
It seems that
a choice of a point of $Spec({\cal Z}(\Mg ))$ corresponds to a choice of a     
polarization on the moduli space.

It seems that the results obtained in the paper can be generalized to the graph
algebras corresponding to arbitrary quantized universal enveloping algebras.

{\bf Acknowledgements:} The author would like to thank  
G.Arutyunov, P.Schupp and  A.A.Slavnov for discussions. He is  
grateful to Professor J.Wess for kind hospitality and the Alexander von 
Humboldt Foundation for the support.  This work has been supported in part 
by the Russian Basic Research Fund under grant 
number 94-01-00300a.

\section*{Appendix}

Let matrices $A$ and $B$ satisfy the  commutation relations of $\Lo$
\[ 
A^{1}\Rp A^{2}\Rpi=\Ra A^{2}\Rai A^{1},\quad 
B^{1}\Rp B^{2}\Rpi=\Ra B^{2}\Rai B^{1} \]
\[ A^{1}\Rp B^{2}\Rpi=\Rp B^{2}\Rai A^{1} \]
\[ \detq A =a_{11}a_{22}-q^2 a_{21}a_{12} =\l_a ,\quad 
\detq B =b_{11}b_{22}-q^2 b_{21}b_{12} =\l_b \]
These relations can be rewritten in components as follows
\bea
&& a_{11}a_{12}=q^{-2}a_{12}a_{11},\quad a_{11}a_{21}=q^{2}a_{21}a_{11},\quad 
a_{11}a_{22}=a_{22}a_{11}\nonumber\\
&& [a_{12},a_{21}]=-(1-q^{-2})a_{11}(a_{11}-a_{22})\Leftrightarrow 
 a_{12}a_{21}=q^{2}a_{21}a_{12}+(1-q^{-2})(\l_a - a_{11}^2 )\nonumber\\
&& [a_{12},a_{22}]=-(1-q^{-2})a_{11}a_{12},\quad 
[a_{21},a_{22}]=(1-q^{-2})a_{21}a_{11}
\nonumber
\eea
and the same relations for $b_{ij}$.
\bea
&& a_{11}b_{11}=qb_{11}a_{11},\quad a_{11}b_{12}=q^{-1}b_{12}a_{11},\quad 
a_{11}b_{21}=qb_{21}a_{11}+ (q-q^{-1})b_{11}a_{21}\nonumber\\
&& a_{11}b_{22}=q^{-1}b_{22}a_{11}+q^{-1}(q-q^{-1})^2 b_{11}a_{11}+
(q-q^{-1})b_{12}a_{21}
\nonumber
\eea
\bea
&& a_{12}b_{11}=qb_{11}a_{12}+(q-q^{-1})b_{12}a_{11},\quad 
a_{12}b_{12}=qb_{12}a_{12}\nonumber\\
&& a_{12}b_{21}=q^{-1}b_{21}a_{12}+q^{-1}(q-q^{-1})^2 b_{12}a_{21}\nonumber\\
&&\qquad\qquad 
+q^{-2}(q-q^{-1})(b_{22}a_{11}+b_{11}a_{22}+(q^{-2}-2)b_{11}a_{11})\nonumber\\
&& a_{12}b_{22}=q^{-1}b_{22}a_{12}+q^{-1}(q-q^{-1})^2 b_{11}a_{12}+
(q-q^{-1})b_{12}a_{22}-q^{-2}(q-q^{-1})b_{12}a_{11}
\nonumber
\eea
\bea
&&a_{21}b_{11}=q^{-1}b_{11}a_{21},\quad 
a_{21}b_{12}=q^{-1}b_{12}a_{21}+q^{-2}(q-q^{-1})b_{11}a_{11},\quad 
a_{21}b_{21}=qb_{21}a_{21}\nonumber\\
&& a_{21}b_{22}=qb_{22}a_{21}+(q-q^{-1}) b_{21}a_{11}
\nonumber
\eea
\bea
&& a_{22}b_{11}=q^{-1}b_{11}a_{22}+q^{-1}(q-q^{-1})^2 b_{11}a_{11}+
(q-q^{-1})b_{12}a_{21}\nonumber\\
&& a_{22}b_{22}=qb_{22}a_{22}-q^{-3}(q-q^{-1})^2 b_{11}a_{11}+
(q-q^{-1})b_{21}a_{12}-q^{-2}(q-q^{-1})b_{12}a_{21}\nonumber\\
&& a_{22}b_{21}=q^{-1}b_{21}a_{22}+q^{-1}(q-q^{-1})^2 b_{21}a_{11}+
(q-q^{-1})b_{22}a_{21}-q^{-2}(q-q^{-1})b_{11}a_{21}\nonumber\\
&& a_{22}b_{12}=qb_{12}a_{22}+(q-q^{-1})b_{11}a_{12}
\label{abcom}
\eea
Let matrices $C$ and $D$ satisfy the following relations
\[
C^{1}\Rp D^{2}\Rpi=\Rp D^{2}\Rpi C^{1}
\]
The relations look in components as follows
\bea
&& c_{11}d_{11}=d_{11}c_{11}-q(q-q^{-1})d_{12}c_{21},\quad 
c_{11}d_{21}=d_{21}c_{11}+ q(q-q^{-1})(d_{11}-d_{22})c_{21}\nonumber\\
&& c_{11}d_{12}=d_{12}c_{11},\quad 
c_{11}d_{22}=d_{22}c_{11}+ q^{-1}(q-q^{-1})d_{12}c_{21}
\nonumber
\eea
\bea
&& c_{12}d_{11}=d_{11}c_{12}+q(q-q^{-1})d_{12}(c_{11}-c_{22}),\quad 
c_{12}d_{12}=q^2 d_{12}c_{12}\nonumber\\
&& c_{12}d_{21}+q^{-1}(q-q^{-1}) c_{11}(d_{11}-d_{22}) =
q^{-2}d_{21}c_{12}+q^{-1}(q-q^{-1}) (d_{11}-d_{22})c_{22}\nonumber\\
&& c_{12}d_{22}=d_{22}c_{12}-q^{-1}(q-q^{-1})d_{12}(c_{11}-c_{22})
\nonumber
\eea
\[
c_{21}d_{11}=d_{11}c_{21},\quad 
c_{21}d_{12}=q^{-2}d_{12}c_{21},\quad 
c_{21}d_{21}=q^2 d_{21}c_{21},\quad 
c_{21}d_{22}=d_{22}c_{21}
\]
\bea
&& c_{22}d_{11}=d_{11}c_{22}+q^{-1}(q-q^{-1})d_{12}c_{21},\quad 
 c_{22}d_{22}=d_{22}c_{22}-q^{-3}(q-q^{-1})d_{12}c_{21}\nonumber\\
&& c_{22}d_{12}=d_{12}c_{22},\quad 
c_{22}d_{21}=d_{21}c_{22}-q^{-1}(q-q^{-1})(d_{11}-d_{22})c_{21}
\label{cdcom}
\eea
The nontrivial commutation relations of the elements $c_{ij}$ with 
$d_{11}^{\frac {1} {2}}$, which are used to define $\Lg^*$, are given by the
formulas
\[ c_{11}d_{11}^{\frac {1} {2}}=d_{11}^{\frac {1} {2}}c_{11}+(1-q)
d_{11}^{-\frac {1} {2}}d_{12}c_{21}\]
\[ c_{12}d_{11}^{\frac {1} {2}}=d_{11}^{\frac {1} {2}}c_{12}-(1-q)
d_{11}^{-\frac {1} {2}}d_{12}(c_{11}-c_{22})+q^{-3}(1-q)^2
d_{11}^{-\frac {3} {2}}d_{12}^2c_{21}\]

Let matrices $A$ and $M$ have the commutation relations
\[
M^{1}\Rp A^{2}\Rpi=\Ra A^{2}\Rai M^{1}
\]
The relations look in components as follows
\bea
&& a_{11}m_{11}=m_{11}a_{11},\quad 
a_{11}m_{12}=m_{12}a_{11}- q(q-q^{-1})m_{11}a_{12}\nonumber\\
&&a_{11}m_{21}=m_{21}a_{11}+ q^{-1}(q-q^{-1})m_{11}a_{21}\nonumber\\
&&a_{11}m_{22}=m_{22}a_{11}+q(q-q^{-1})(m_{12}a_{21}-m_{21}a_{12})- 
(q-q^{-1})^2 m_{11}(a_{22}-a_{11})
\nonumber
\eea
\bea
&& a_{12}m_{11}=q^2 m_{11}a_{12},\quad 
a_{12}m_{12}=m_{12}a_{12}\nonumber\\ 
&&a_{12}m_{21} =
m_{21}a_{12}-q^{-1}(q-q^{-1})m_{11}(a_{11}-a_{22})\nonumber\\
&& a_{12}m_{22}=q^{-2}m_{22}a_{12}-q^{-1}(q-q^{-1})m_{12}(a_{11}-a_{22})+
q^{-1}(q-q^{-1})(q^2 -q^{-2})m_{11}a_{12}
\nonumber
\eea
\bea
&&a_{21}m_{11}=q^{-2}m_{11}a_{21},\quad 
a_{21}m_{12}=m_{12}a_{21}+q^{-1}(q-q^{-1})m_{11}(a_{11}-a_{22})\nonumber\\
&&a_{21}m_{21}=m_{21}a_{21},\quad 
a_{21}m_{22}=q^2 m_{22}a_{21}+q(q-q^{-1})m_{21}(a_{11}-a_{22})
\nonumber
\eea
\bea
&& a_{22}m_{11}=m_{11}a_{22},\quad a_{22}m_{12}=m_{12}a_{22}+
q^{-1}(q-q^{-1})m_{11}a_{12}\nonumber\\
&&a_{22}m_{21}=m_{21}a_{22}-q^{-3}(q-q^{-1})m_{11}a_{21}\nonumber\\
&&
a_{22}m_{22}=m_{22}a_{22}-q^{-1}(q-q^{-1})(m_{12}a_{21}-m_{21}a_{12})\nonumber\\
&&\qquad\qquad +q^{-2}(q-q^{-1})^2 m_{11}(a_{22}-a_{11})
\label{amcom}
\eea
The nontrivial commutation relations of the elements $m_{ij}$ with 
$a_{11}^{\frac {1} {2}}$, which are used to define $\Lg^*$, are given by the
formulas
\[ a_{11}^{\frac {1} {2}}m_{12}=m_{12}a_{11}^{\frac {1} {2}}
+(1-q)m_{11} a_{11}^{-\frac {1} {2}}a_{12}\]
\[ a_{11}^{\frac {1} {2}}m_{21}=m_{21}a_{11}^{\frac {1} {2}}
+(1-q^{-1})m_{11} a_{11}^{-\frac {1} {2}}a_{21}\]

\end{document}